# Cooperative Autonomous Vehicle Speed Optimization near Signalized Intersections

Mahmoud Faraj, Baris Fidan, and Vincent Gaudet
Faculty of Engineering, University of Waterloo, ON, Canada
Email: {msfaraj,fidan,vcgaudet}@uwaterloo.ca

## Abstract

Road congestion in urban environments, especially near signalized intersections, has been a major cause of significant fuel and time waste. Various solutions have been proposed to solve the problem of increasing idling times and number of stops of vehicles at signalized intersections, ranging from infrastructure to vehicle-based techniques. However, all the solutions introduced to solve the problem have approached the problem from a single vehicle point of view. This research introduces a game-theoretic cooperative speed optimization framework to minimize vehicles' idling times and number of stops at signalized intersections. This framework consists of three modules to cover individual autonomous vehicle speed optimization; conflict recognition; and cooperative speed decision making. A time token allocation algorithm is introduced through the proposed framework to allow the vehicles to cooperate and agree on certain speed actions such that the average idling times and number of stops are minimized. Simulation to test and validate the proposed framework is conducted and results are reported.

## 1 Introduction

In the U.S. every year, 4.8 billion hours are wasted by traffic congestion [1]. Some automobile industrial companies have made leaps toward manufacturing Autonomous Vehicles (AVs), which can implement different levels of automatic functions as a means of achieving a safer and more efficient transportation system [2–6]. Some research has been conducted on optimal speed computation near Signalized Intersections (SIs) to minimize the vehicles' idling times and number of stops. [7] has proposed an algorithm called Green Light Optimal Speed Advisory (GLOSA) to minimize the number of stop times at SIs through a journey. The impacts of this algorithm on traffic efficiency and average trip time were reported. The performance analysis of the same algorithm was investigated again in [8] using the performance metrics of average fuel consumption and average stop times at SIs.

[9] has introduced an approach for a single vehicle optimal speed computation to

reduce fuel consumption and $CO_2$ emissions. In this approach though, the case of a driver cruising his/her vehicle speed to pass a TL is not considered. [10] has used the approach proposed in [9] to reduce $CO_2$ emissions. Since [10] pays attention to the reduction of $CO_2$ by reducing stop-and-go driving, the case where the driver may cruise his/her vehicle speed to pass the TL is taken into account. [11] has investigated the impacts of Vehicle-to-Infrastructure (V2I) communication, namely TL to vehicle communication, on fuel and emission reductions. [12] has posed a speed advisory model that computes a fuel-optimal speed profile during deceleration and acceleration phases at SIs.

Game theory was introduced in the early years of the twentieth century in [13] and [14]. It has been applied to some transportation problems, such as shortest path and road congestion problems. [15] has proposed a shortest-path game with transferable utility, focusing on the allocation of profits generated by the coalitions of players. [16] has presented a shortest-path game in which players own road segments in a network. Each player in the game receives a non-negative reward if he/she transports a good from the source to the destination. [17] has introduced a model in which players who share resources (i.e., routes) can form coalitions to selfishly compete against each other to maximize their values. [18] has discussed the similarities between cooperative congestion games and their non-cooperative counterparts to demonstrate important issues, such as the existence of and the convergence to a pure strategy Nash Equilibrium (NE).

Furthermore, game theory has been applied to the dynamic TL signal timing control problem. [19] has introduced a model for TL system control based on Markov Chain game with the objective of minimizing the queue lengths at multiple SIs. [20] has proposed a two-player cooperation game for TL signal timing control applied to a two-phase SI. Similar research to [19] is presented in [21] where a non-cooperative game to model the TL signal timing control problem is introduced based on game theory and modeled as a finite controlled Markov Chain. However, the TL model in [21] is applied to a single SI. [22] has presented a game theory model based on Cournot's Oligopoly game. [23] has proposed a novel game theory optimization algorithm for TL signal timing control. The Nash Bargaining (NB) is used to find the optimal strategy of the TL signal timing.

Noticeably, all the techniques proposed for minimizing the idling times and number of stops of vehicles approaching an SI have been limited to optimal speed computation in which the vehicles individually and independently compute and update their optimal speeds. In this research work, a Cooperative Speed Optimization Framework (CSOF) is proposed. The proposed CSOF is designed to function on AVs of the highest level of autonomy (i.e., the AV is completely autonomous and the driver does not take control of the AV at any point in time). The cooperative notion of the CSOF is to allow the AVs to interact with each other and agree on implementing certain speed actions when approaching SIs. It relies on linear programming and game theory, consisting of three modules to cover individual vehicle speed optimization, conflict recognition, and speed optimization decision making. A time token allocation algorithm is proposed to be embedded in the TL such that AVs are able to cooperate with each other and with the TL. Thus, the average idling times and number of stops at SIs are minimized.

The rest of the paper is organized as follows: Section 2 states the problem of minimizing the idling times and number of stops of AVs approaching an SI as well as addresses the game theory formulation toward the proposed solution. Section 3 introduces the game theory cooperative framework, CSOF, proposed to minimize the average idling times and number of stops of AVs. Section 4 poses the concept of cooperative bargaining. Section 5 presents the simulation environment and reports results investigating the performance of the CSOF. Finally, concluding remarks and future work recommendations are discussed in Section 6.

## 2 Problem Statement and Game-theoretic Formulation

### 2.1 Signalized Roadway Intersection Setting

Consider as a simple example the illustration of a two-lane, four-roadway SI (Fig. 1). For simplicity, assume that the TL control system has a two-phase static cycle where East and West roadways are one phase and North and South roadways are the other phase. Each phase has a signal design of Green-Yellow-Red; however, for simplicity, the yellow-light time is assumed to be part of the green-light-time duration. The key TL parameters are the green-light-time duration $T_g$, the red-light-time duration $T_r$, and the TL cycle duration $T_c = T_g + T_r$. These parameters are assumed to be constant, e.g., $T_g = 24sec$, $T_r = 36sec$, and $T_c = 60sec$. Assume that there is V2I communication such that the vehicles heading toward the TL can receive signal timing information and that every AV is conducting *speed optimization re-planning* to have a chance of meeting the green-light time.

**Definition 1** *Speed optimization re-planning is a game in which each AV performs speed optimization every time step t. There is a probability p that the AV will proceed according to the previous strategy at time step $t-1$ and a probability $(1-p)$ that it will move to a different strategy (i.e., adopt a different speed).*

To clarify the complexity of the problem, assume that for a certain cycle, the arrival and maximum departure rates of each roadway at the TL are $\lambda = 0.25veh/sec$ and $\mu = 0.333veh/sec$, respectively [24]. Therefore, in this particular cycle, the number of AVs arriving from each roadway during the red time is $N_{arr} = \lambda T_r = (0.25)(36) = 9veh$, while the maximum number of AVs that can depart the TL from each roadway during the green time and whole cycle is $N_{dep} = \mu T_g = (0.333)(24) = 8veh$.

Making this setting, assume that the TL has just turned green for the East-West directions. Consider the case of two AVs travelling on the West roadway performing *speed optimization re-planning*. Taking the queue size into account, according to the computations by these AVs, each of them can pass within the current green light. Since only one AV can pass through, the other will experience an unexpected delay, waiting for the next green light. Hence, AVs negatively impact the objectives of each other.

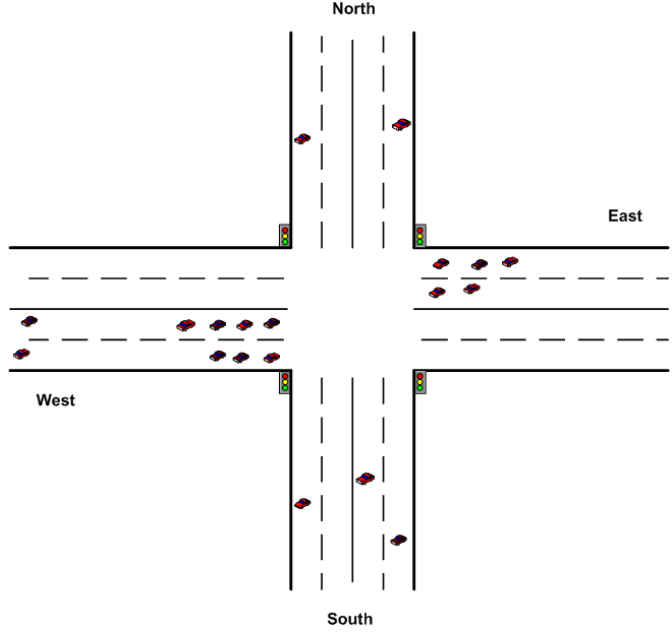

**Figure 1:** A simple example of a traffic light scenario.

## 2.2 Problem Formulation

Consider a group of AVs travelling within a locality with $m$ TLs. Each AV with index $AV_i$ contemplates speed optimization to minimize its idling times at TLs from an initial location $p_{0i}$ to a final destination $p_{fi}$. The trip from $p_{0i}$ to $p_{fi}$ is made on a path, $P(p_{0i}, p_{fi})$, constructed from a set of road segments ending with TLs, $L = \{L_1, L_2, \ldots, L_m\}$. The speed $v_i(t)$ of each $AV_i$ belongs to a set of feasible speeds, $V = \{\bar{v}_1, \bar{v}_2, \ldots, \bar{v}_f\}$. The cost of the trip for $AV_i$ on a road segment $L_j$, where $j = 1, 2, \ldots, m$, denoted by $C_{sv}^{L_j}(i)$, explicitly models the idling times of AVs. For $AV_i$, the cost of travelling over road segment $L_j$ incurred by choosing a time indexed sequence of velocities, $sv$, is defined as follows:

$$C_{sv}^{L_j}(i) = \begin{cases} t_i & if \quad stop \\ 0 & if \quad no\, stop \end{cases} \quad (1)$$

where $t_i$ is the idling time of $AV_i$ at the TL positioned at the end of road segment $L_j$. The total cost for $AV_i$, incurred over a path $P(p_{0i}, p_{fi})$ that is composed of the road segments $L_1^P, \ldots, L_{N_P}^P \in L$ (sequentially), is the summation of idling times at all TLs along the path, i.e.,

$$C_{sv}^{P}(i) = \sum_{j=1}^{N_P} C_{sv}^{L_j^P}(i). \tag{2}$$

where $N_P$ is the number of TLs on the path $P$, and $sv$ denotes the sequence of velocities for road segment $L_j^P$. $sv$ is the concatenation of $sv_1, \ldots, sv_{N_P}$. To provide a sub-optimal solution to the above overall task, we follow a decentralized approach and consider each road segment $L_j \in L$ of the locality separately. We propose a "Time Token Allocation Algorithm" for the TL and a cooperative distributed conflict resolution scheme for the vehicles in each such $L_j$.

For player $AV_i$, the optimal speed value in the set of possible speeds may lead to a time token within the green light, allocated by the TL, $\tau_i$. The time token $\tau_i$ is the index of a time window assigned by the TL using the Time Token Allocation Algorithm. This means that for player $AV_i$, the cost associated with the time token $\tau_i$ is the minimum (e.g., player $AV_i$ will pass through the TL without stopping, $C_{sv}^{L_j}(i) = 0$).

We define a cooperative speed optimization game, $G$. In this game, AVs with conflicting allocated time tokens agree to take certain speed actions to resolve the conflict. Thus, for each player $AV_i$, there is a finite non-empty set of speed actions $V$. There is an idling time cost, $C_{sv}^{L_j}(i)$, associated with each sequence of actions $sv$. Action sequences are associated with a preference relationship such that $C_{sv^*}^{L_j}(i) < C_{sv}^{L_j}(i)$ means $sv^* \succ sv_j$ (i.e., the sequence of actions $sv^*$ is preferred over that $sv$ as it incurs less cost). As such, the speed optimization game is represented as $G = (n, V, C)$.

## 2.3 Stability of the Speed Optimization Game

In the *speed optimization re-planning* game (Definition 1), the chosen actions of AVs are considered pure strategies, and therefore, there always exists pure equilibrium where no AV wishes to unilaterally change its optimal solution. However, this is true only for non-strictly competitive games. As the game progresses, AVs compete to gain resources such that one AV's gain is another AV's loss.

A mixed strategy game, which always has a mixed equilibrium, is a game in which the strategies available to the players are not deterministic but are regulated by probabilistic rules [25]. Thus, from Definition 1, it is concluded that there is a probability distribution over all the strategies available to every AV in the game. Hence, the *speed optimization re-planning* game, as described in Definition 1, is a mixed strategy game for which a mixed equilibrium always exists.

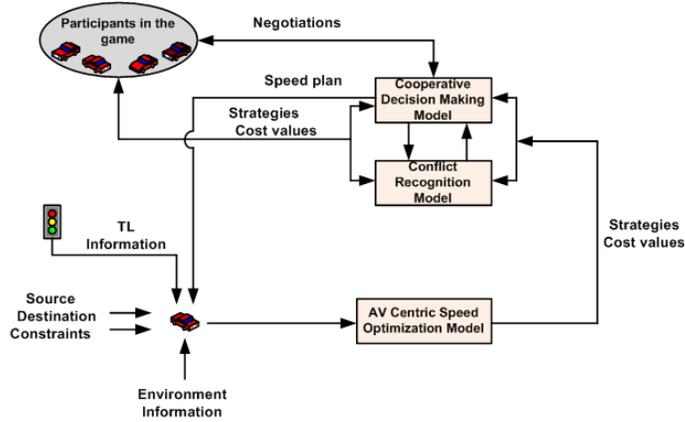

**Figure 2:** Schematic depiction of the cooperative speed optimization framework.

# 3 The Proposed Cooperative Speed Optimization Framework

Schematics of the Cooperative Speed Optimization Framework (CSOF) we propose to solve the *speed optimization re-planning* game is provided in Fig. 2. This framework consists of three modules to cover issues of (i) AV rational speed optimization, (ii) information and conflict recognition, and (iii) cooperative speed optimization decision making. However, before proceeding to the core of the proposed scheme, we address a few issues with respect to the safety constraints of consecutive AVs on the roadway.

## 3.1 Car Following Model

The CSOF is designed to function on multiple-lane roadways only. A car-following model is defined with two essential rules. First, AVs using the CSOF in free motion can smoothly overtake each other on the roadway to comply to certain speed actions resulting from their interaction and cooperation. Second, under certain traffic conditions such as when overtaking is not possible, a safe following distance between consecutive AVs is maintained. Consecutive AVs are modeled to maintain a minimum time gap of two seconds in order to avoid collision [26]. All the AVs are identical in length and have an average length of $5 meters$. The reaction time of AVs to the sudden deceleration of the traffic ahead is assumed to be $1.10 seconds$ [27].

## 3.2 Autonomous Vehicle Speed Optimization Module

The objective of this module is to provide each player (vehicle), $AV_i$, with the optimal speed at every time step $t$. Based on the Time to Intersection $TTI_i$ and using a Time Token Allocation Algorithm, the TL may allocate $\tau_i$ to $AV_i$, where $\tau_i$ is an integer value indicating the index of a time window during which $AV_i$ can pass the intersection smoothly. The speed $v_i(t)$ of $AV_i$ at time step $t$ is a function of the traffic density $D(L_j)$ on road segment $L_j$. [28] justified the linearity of the relation between traffic density and speed under mild generic assumptions, concluding that as traffic concentration/density increases, speed decreases. As such, some definitions are stated.

- The maximum speed $AV_i$ can travel at, $v_{max}$, will only occur when there are no other vehicles on the roadway.

- In general, the speed of $AV_i$ goes to zero as the road reaches the maximum density, $v_i(t)$ converges to 0 as $D(L_j)$ converges to $D_{max}(L_j)$.

Therefore, considering a linear relation between the traffic density and speed, the speed $v_i(t)$ of $AV_i$ with respect to the traffic density $D(L_j)$ on road segment $L_j$ is

$$v_i(t) = v_{max} \left(1 - \frac{D(L_j)}{D_{max}(L_j)}\right) \tag{3}$$

$AV_i$ is allocated a token $\tau_i$ only if $TTI_i$ falls within the upcoming green-light time, i.e., $TTI_i \leq R_g$ or $R_r < TTI_i \leq R_r + T_g$ where $R_g$ and $R_r$ are the remaining green-light and red-light times respectively. For $AV_i$ approaching a TL, the speed that minimizes $AV_i's$ idling time cost is found as follows:

### 3.2.1 Light is Green

As $AV_i$ receives upcoming signal information from the TL, indicating that the current light is green, there are three possible cases in terms of $TTI_i$ and $R_g$:

- Case 1: $TTI_i \leq R_g$. In this case, using the current speed, $AV_i$ will be able to pass through within the remaining green-light time. The TL allocates a time token $\tau_i$ to $AV_i$. Thereby, $AV_i$ maintains its speed to pass during the assigned time token.

$$\begin{aligned} s_i &= v_i(t) \\ subject\ to&: \\ v_i(t) &\geq d_i(t)/a_i \\ v_i(t) &\leq d_i(t)/b_i \\ v_i(t) &\geq v_{min} \\ v_i(t) &\leq v_{max} \end{aligned}$$

where $TTI_i = d_i(t)/v_i(t)$, $v_i(t)$ is the speed of $AV_i$ at time step $t$, $d_i(t)$ is the distance of $AV_i$ to the stop line of the TL at time step $t$, $s_i$ is the optimal speed of $AV_i$ at time step $t+1$, and $v_{max}$ and $v_{min}$ are the maximum and minimum speed limits on the road segment respectively, while $a_i$ and $b_i$ represent the lower and upper boundaries of the allocated time token respectively.

$$a_i = (\tau_i - 1)\frac{1}{\mu} \tag{4}$$

$$b_i = \tau_i \frac{1}{\mu} \tag{5}$$

where $\mu$ is the departure rate in $veh/sec$.

- Case 2: $R_g + T_r \geq TTI_i > R_g$. In this case, the vehicle is not allocated a time token, and the speed of the vehicle is optimized over the distance to the TL so that $TTI_i$ is sufficient to meet the next green light.

$$\begin{aligned}
s_i = min \quad & v_i(t) \\
subject\ to: & \\
& v_i(t) \geq d_i(t)/(R_g + T_r + T_q) \\
& v_i(t) \leq d_i(t)/(R_g + T_r + T_g) \\
& v_i(t) \geq v_{min} \\
& v_i(t) \leq v_{max}
\end{aligned}$$

where $T_q$ is the time needed to clear all the vehicles in the queue, and it is found as follows:

$$T_q = \frac{n(t)}{\mu} \tag{6}$$

where $n(t)$ denotes the number of vehicles currently in the queue. In addition, if the current speed does not allow $AV_i$ to be part of the green-light time but the maximum speed of the roadway does, the speed optimization system will accelerate the speed of $AV_i$ such that it is allocated a token.

$$\begin{aligned}
s_i = max \quad & v_i(t) \\
subject\ to: & \\
& v_i(t) \geq d_i(t)/a_i \\
& v_i(t) \leq d_i(t)/b_i \\
& v_i(t) \geq v_{min} \\
& v_i(t) \leq v_{max}
\end{aligned}$$

- Case 3: $R_g + T_r + T_g \geq TTI_i > R_g + T_r$. In this case, $AV_i$ will maintain its current speed as $TTI_i$ leads $AV_i$ to be part of the green-light time of the next cycle; However, $AV_i$ will not yet be allocated a time token.

$$\begin{aligned}
s_i = & v_i(t) \\
subject\ to: & \\
& v_i(t) \geq v_{min} \\
& v_i(t) \leq v_{max}
\end{aligned}$$

### 3.2.2 Light is Red

If the information received by the vehicle from the TL indicates that the current light is red, there are three possible cases in terms of $TTI_i$ and $R_r$:

- Case 1: $TTI_i < R_r$. In this case, $AV_i$ will not be allocated a time token, and its speed is optimized such that it will meet the next green light.

$$s_i = min \quad v_i(t)$$
$$subject\ to:$$
$$v_i(t) \geq d_i(t)/(R_r + T_q)$$
$$v_i(t) \leq d_i(t)/(R_r + T_g)$$
$$v_i(t) \geq v_{min}$$
$$v_i(t) \leq v_{max}$$

- Case 2: $R_r < TTI_i \leq R_r + T_g$. In this case, $AV_i$ is allocated a time token within the upcoming green-light time.

$$s_i = v_i(t)$$
$$subject\ to:$$
$$v_i(t) \geq d_i(t)/(R_r + a_i)$$
$$v_i(t) \leq d_i(t)/(R_r + b_i)$$
$$v_i(t) \geq v_{min}$$
$$v_i(t) \leq v_{max}$$

- Case 3: $TTI_i > R_r + T_g$. In this case, $AV_i$ is not allocated a time token and its speed is optimized to meet the green-light time of the next cycle.

$$s_i = min \quad v_i(t)$$
$$subject\ to:$$
$$v_i(t) \geq d_i(t)/(R_r + T_g + T_r + T_q)$$
$$v_i(t) \leq d_i(t)/(R_r + T_r + 2T_g)$$
$$v_i(t) \geq v_{min}$$
$$v_i(t) \leq v_{max}$$

### 3.2.3 Energy Consumption Model

All the AVs involved in the cooperative process are assumed to be Electric Autonomous Vehicles (EAVs); therefore, the energy consumption model presented in [29] has been modified to compute the instant energy consumed by every AV at time step $t$. The total energy cost consumed by $AV_i$ at time step $t$ consists of multiple sub-costs as follows:

- Potential Consumed/Gained Energy: the potential energy $EC_i^P(t)$ at time step $t$ is consumed from the battery during the uphill travel and is gained into the battery during the downhill travel. The potential consumed and gained energies are

$$EC_i^{PC}(t) = \frac{1}{\eta}[m\,g\,u(t)] \tag{7}$$

$$EC_i^{PG}(t) = -\frac{1}{\eta}[m\,g\,u(t)] \tag{8}$$

where $\eta$ is the efficiency of $AV_i$, $m$ is the mass of $AV_i$, $g$ is the gravity factor, and $u(t)$ is the elevation of the road segment at time step $t$.

- Loss of Energy: the loss of energy at time step $t$, which is always consumed from the battery, occurs due to aerodynamic and rolling resistances.

$$EC_i^{loss}(t) = \frac{1}{\eta}[f_r\,m\,g\,v_i(t) + \frac{1}{2}\rho\,A\,d_r\,v_i^3(t)] \tag{9}$$

where $f_r$ is the friction coefficient, $\rho$ is the air density coefficient, $A$ is the cross sectional area of $AV_i$, and $d_r$ is the air drag coefficient.

- Acceleration/Deceleration Energy: the acceleration energy $EC_i^{ac}(t)$ at time step $t$ is consumed from the battery as $AV_i$ accelerates to a higher speed while the deceleration energy $EC_i^{dc}(t)$ at time step $t$ is recuperated and stored into the battery as $AV_i$ comes to a lower speed.

$$EC_i^{ac}(t) = \frac{1}{\eta}P_i^{wr}\frac{d_i(t)}{v_i(t) - v_i(t-1)} \tag{10}$$

$$EC_i^{dc}(t) = -\frac{1}{\eta}P_i^{wr}\frac{d_i(t)}{v_i(t) - v_i(t-1)} \tag{11}$$

where $P_i^{wr}$ is the power of the electric motor of $AV_i$.

- Energy Consumed by On-Board Electric Devices: this energy is not path related and is consumed directly from the battery at time step $t$ by the on-board electric devices such as air conditioner, windshield wipers, etc.

$$EC_i^{ed}(t) = \sum_{j=1}^{n} P_i^{ed}(j)\,t_i^{ed}(j) \tag{12}$$

where $P_i^{ed}(j)$ is the power withdrawn at time step $t$ by the electric device $j$ and $t_i^{ed}(j)$ is the time that device $j$ takes in use.

Therefore, the total energy cost consumed by $AV_i$ at time step $t$ is computed as

$$\begin{aligned}EC_i^T(t) = [&EC_i^{PC}(t) + EC_i^{PG}(t) + EC_i^{loss}(t)\\&+ EC_i^{ac}(t) + EC_i^{dc}(t)]\end{aligned} \tag{13}$$

## 3.3 Information and Conflict Recognition Module

In this module, if two or more players have been allocated the same time token, the TL informs them that they have a conflict. Players with conflicting tokens communicate with each other to share their strategies and associated costs. Consequently, they start to negotiate to find a binding agreement based on which they can cooperate and agree on certain speed actions. Once an agreement is reached, all the players abide by the rules to apply those actions.

## 3.4 Cooperative Decision Making Module

The cooperative game notion in this module is based on the assumption that players can reach a binding agreement with which they commit to apply certain strategic actions. As players are assumed to be rational, if the idling time cost to a player is greater than what it would have been without cooperation, then the rationality assumption is violated. Therefore, the following rationality axiom is stated.

**Axiom 1** *At time step t, there exists an optimal strategy $v_k$ for player $AV_i$ such that $C_{v_k}^{L_j}(i) \leq C_{v_j}^{L_j}(i), \forall k, j \in V$ (i.e., the cost associated with this optimal strategy is less than or equal to that associated with any other strategy player $AV_i$ can take). However, player $AV_i$ is free to choose any other strategy that might yield a higher cost, but only in exchange for a reward.*

The TL allocates time tokens to the players using Algorithm 1. In Algorithm 1, *VIN* is the AV identification number, $N_q$ is the number of vehicles currently in the queue, *Tsd* is the slot duration (i.e., time token duration), and *Tslot* is the time token location in the TL memory. According to this algorithm, the TL gives priority in allocating tokens to the queued AVs. The rest of the green light time is segmented as tokens and offered to the approaching AVs. Players with conflicting tokens will list the costs caused by the conflict rather than the expected ones. It is assumed that each AV has a mode property, which may take one of three values at a time: *Rush Mode*, *Normal Mode*, or *Relaxed Mode*.

- *Rush Mode*: this mode is used by the AV for urgent and emergency situations (e.g., must be in the hospital shortly).

- *Normal Mode*: this mode is used by the AV when there is no emergency; the AV may yield the road to others.

- *Relaxed Mode*: this mode is used by the AV when there is plenty of time. The AV would yield the road to other vehicles comfortably.

Each mode is represented by a scalar value. For instance, the mode values may be 0, 1, and 2 for modes *Relaxed*, *Normal*, and *Rush* respectively. Players will first play the

**Algorithm 1** : Time Token Allocation Algorithm for a Single Roadway

*Input:* $VIN_i, TTI_i$ ; *Output:* $\tau_i$

  **if** Light is Green **then**
    **if** $(TTI_i \leq R_g)$ **then**
      **for** $j = N_q + 1 : N_{dep}$ **do**
        $a = Tsd * (j-1)$;
        $b = Tsd * j$;
        **if** $(TTI_i \geq a\ \&\ TTI_i \leq b)$ **then**
          $Tslot(1, j) = VIN_i$;
          $\tau_i = j$;
          *break*;
        **end if**
      **end for**
    **else**
      $\tau_i = 0$;
    **end if**
  **else**
    **if** $(TTI_i < R_r)$ **then**
      $\tau_i = 0$;
    **else if** $(TTI_i > R_r + (Tsd * N_q)\ \&\ TTI_i \leq R_r + T_g)$ **then**
      **for** $j = N_q + 1 : N_{dep}$ **do**
        $a = R_r + Tsd * (j-1)$;
        $b = R_r + Tsd * j$;
        **if** $(TTI_i \geq a\ \&\&\ TTI_i \leq b)$ **then**
          $Tslot(1, j) = VIN_i$;
          $\tau_i = j$;
          *break*;
        **end if**
      **end for**
    **else**
      $\tau_i = 0$;
    **end if**
  **end if**

game based on the mode type. The one using the mode with the highest value will win (i.e., use the current time token), while the one using the mode with a smaller value will lose (i.e., slow down its own speed and request a different time token). However, the TL will grant the loser a credit point and deduct a credit point from the winner. The winner of the mode-based game is determined using the following formula.

$$AV_{winner} = max(M(AV_1), M(AV_2)) \quad (14)$$

If both players are using the same mode, they will decide the winner based on the credit points that they have. The one with the most will eventually win the game. Again, a credit point is deducted from the winner, and a credit point is granted to the loser. In this case, the winner is determined using the following formula.

$$AV_{winner} = max(CP(AV_1), CP(AV_2)) \quad (15)$$

If it happens that both players have the same mode value and number of credit points, a random number-generation procedure between the TL and players is conducted to resolve the issue. Basically, each of the players as well as the TL will generate a random number. The one whose generated number is closer to that of the TL will win the current time token but lose a credit point. The other will gain a credit point and request a different token. The winner of the game is determined using the following formula.

$$RN_1 = |RN(TL) - RN(AV_1)|$$
$$RN_2 = |RN(TL) - RN(AV_2)|$$

$$AV_{winner} = \begin{cases} AV_1 & if \quad RN_1 < RN_2 \\ AV_2 & otherwise \end{cases} \quad (16)$$

To further clarify the cooperative speed optimization game, an example of two AVs approaching a TL is presented next.

**Example 1** *Recall the problem statement example (Fig. 1) and assume that the TL has just turned green for the East-West directions. After communicating with the TL, the two AVs, approaching the TL from the West, have been allocated, at time step t, the same and only-remaining time token. Both vehicles will have only two strategies to choose from. For instance, player $AV_1'$s available strategies are $V^{AV_1} = \{v_1(t)^{AV_1}, v_2(t)^{AV_1}\}$, where using strategy $v_1(t)^{AV_1}$ corresponds to using the current time token and using strategy $v_2(t)^{AV_1}$ corresponds to minimizing its own speed and requesting a token within the next green light. Table 1 has been constructed to clarify the game setting.*

*When both players choose the same strategy resulting in the use of strategy profile $(C^{L_j}_{v_1(t)}(1), C^{L_j}_{v_1(t)}(2))$ or $(C^{L_j}_{v_2(t)}(1), C^{L_j}_{v_2(t)}(2))$, the cost is high for both of them. When*

Table 1: Action and response table of a two-player game

| Players/Strategies | $v_1(t)^{AV_2}$ | $v_2(t)^{AV_2}$ |
|---|---|---|
| $v_1(t)^{AV_1}$ | $(C^{L_j}_{v_1(t)}(1), C^{L_j}_{v_1(t)}(2))$ | $(C^{L_j}_{v_1(t)}(1), C^{L_j}_{v_2(t)}(2))$ |
| $v_2(t)^{AV_1}$ | $(C^{L_j}_{v_2(t)}(1), C^{L_j}_{v_1(t)}(2))$ | $(C^{L_j}_{v_2(t)}(1), C^{L_j}_{v_2(t)}(2))$ |

*either of the two players chooses the optimal strategy while the other chooses the second preferred strategy resulting in the use of strategy profile $(C^{L_j}_{v_2(t)}(1), C^{L_j}_{v_1(t)}(2))$ or $(C^{L_j}_{v_1(t)}(1), C^{L_j}_{v_2(t)}(2))$, the game is stable. In this case, the strategy profile is an NE and Pareto Optimal (PO). Hence, the binding agreement between the players would enforce the idea that they should choose different strategies. A numerical representation of such a game may be shown as in Table 2.*

Table 2: Action and response table of a two-player game

| Players/Strategies | $v_1(t)^{AV_2}$ | $v_2(t)^{AV_2}$ |
|---|---|---|
| $v_1(t)^{AV_1}$ | (4,4) | (0,2) |
| $v_2(t)^{AV_1}$ | (2,0) | (3,3) |

### 3.4.1 Multi-Phase Cooperative Speed Optimization Game

When more than two vehicles are allocated the same token, a multi-phase cooperative procedure is implemented to resolve the conflict. The multi-phase game is composed of multiple two-player sub-games. In each sub-game, only two players cooperate to find their acceptable joint strategies. Then, the winners of the two-player sub-games will play another sub-game and so on until the winner of the only available time token is determined.

## 4 Cooperative Credit-Point Bargaining

A cooperative game is a tuple of two elements $(N, f)$, where $N = \{1, 2, ...n\}$, is a finite set of AVs willing to trade credit points, and $f$ is a function that maps subsets of $N$ to numbers. If $S$ is a subset of $N$, such that $S \subseteq N$, then $f(S)$ is the total value induced when the members of $S$ come together to trade credit points. For further clarification, an example is presented next.

**Example 2** *Assume that there are three AVs, $N = \{AV_1, AV_2, AV_3\}$, heading toward a TL, where $AV_1$ is a seller of a credit point, while $AV_2$ and $AV_3$ are two buyers. Consider the case that $AV_1$ has only one credit point to sell at \$3 and each of the buyers contemplates to buy at most one credit point. $AV_2$ is willing to pay \$5, while $AV_3$ is willing to* *pay $8. The characteristic function, f, of this game is defined as follows:*

$$f(AV_1) = f(AV_2) = f(AV_3) = 0$$
$$f(AV_1, AV_2) = 5 - 3 = 2$$
$$f(AV_1, AV_3) = 8 - 3 = 5$$
$$f(AV_2, AV_3) = 0$$
$$f(AV_1, AV_2, AV_3) = 8 - 3 = 5$$

## 4.1 The Marginal Contribution

As introduced by [30], the marginal contribution concept provides the analytical reasoning of bargaining. Let $N \backslash AV_i$ be the subset of $N$ that contains all the AVs except $AV_i$. The marginal contribution of $AV_i$ is $f(N) - f(N \backslash AV_i)$ and denoted by $MC_{AV_i}$. For example, the marginal contributions of the previously defined game are

$$MC_{AV_1} = f(N) - f(N \backslash AV_1) = 5 - 0 = 5$$
$$MC_{AV_2} = f(N) - f(N \backslash AV_2) = 5 - 5 = 0$$
$$MC_{AV_3} = f(N) - f(N \backslash AV_3) = 5 - 2 = 3$$

**Definition 2** *An allocation, $(x_{av_1}, x_{av_2}, ..., x_{av_n})$, which is a collection of numbers representing the division of the overall value, where $x_{av_i}$ indicates the value received by $AV_i$, is individually rational if $x_{av_i} \geq f(AV_i)$, $\forall\, i \in \{1, 2, ..., n\}$.*

**Definition 3** *An allocation, $(x_{av_1}, x_{av_2}, ..., x_{av_n})$, is efficient if $\sum_{i=1}^{n} x_{av_i} = f(N)$*

**Definition 4** *An individually rational and efficient allocation, $(x_{av_1}, x_{av_2}, ..., x_{av_n})$, satisfies the Marginal Contribution Principle if $x_{av_i} \leq MC_{AV_i}$, $\forall\, i \in \{1, 2, ..., n\}$.*

## 4.2 The Core

The core is the solution concept of coalitional games, containing all the possible allocations (i.e., divisions of the overall value) [33]. Let $x(S)$ be the sum of the values received by the AVs in the subset $S$, such that

$$x(S) = \sum_{i \in S} x_{av_i} \qquad (17)$$

In addition, let the marginal contribution of a subset $S$ of $N$ be $MC_S = f(N) - f(N \backslash S)$. According to [30], the core has two main conditions.

**Theorem 1** *An allocation, $(x_{av_1}, x_{av_2}, ..., x_{av_n})$, is part of the core if it is efficient and for every subset S of N, $x(S) \geq f(S)$ is satisfied.*

**Proof 1** *An allocation that belongs to the core of the game is individually rational.*

*let $S = \{AV_i\}$ for $i = 1, 2, ..., n$*

*Noticeably, $x\{AV_i\} = x_{av_i}$ both represent the values received by $AV_i$.*

*Therefore, the condition $x(S) \geq f(S)$ is in fact the individual rationality condition $x_{av_i} \geq f(AV_i)$.*

□

**Theorem 2** *An allocation, $(x_{av_1}, x_{av_2}, ..., x_{av_n})$, is part of the core if it is efficient and for every subset S of N, $x(S) \leq MC_S$ is satisfied.*

**Proof 2** *An allocation that belongs to the core of the game satisfies $x(S) \leq MC_S$.*

*Using the individual rationality condition, consider $N \backslash S$*

$$x(N \backslash S) \geq f(N \backslash S) \tag{18}$$

$$x(N \backslash S) = x(N) - x(S) \tag{19}$$

*By efficiency, we have*

$$x(N) = f(N). \tag{20}$$

*Substituting (19) into (18)*

$$x(N) - x(S) \geq f(N \backslash S) \tag{21}$$

*Substituting (20) into (21)*

$$x(S) \leq f(N) - f(N \backslash S) = MC_S \tag{22}$$

□

Therefore, the core of the cooperative credit point bargaining game is defined as follows:

$$\{(x_{av_1}, x_{av_2}, ..., x_{av_n}) : \sum_{i \in N} x_{av_i} = f(N), \text{ and} \atop x(S) \geq f(S), \forall S \in N\} \tag{23}$$

To find the core elements, we propose that the problem is formulated as a Constraint Satisfaction Problem (CSP). The most common CSP solving techniques are Backtracking Search and Local Search [35]. For instance, the feasible allocations in Example 2 are the points $(x_{av_1}, x_{av_2}, x_{av_3})$, such that

$$x_{av_1} + x_{av_2} + x_{av_3} = 5$$
$$subject\ to:$$
$$x_{av_1} + x_{av_2} \geq 2$$
$$x_{av_1} + x_{av_3} \geq 5$$
$$x_{av_2} + x_{av_3} \geq 0$$
$$x_{av_1} \geq 0,\ x_{av_2} \geq 0,\ x_{av_3} \geq 0$$

The domains of $x_{av_1}$, $x_{av_2}$ and $x_{av_3}$ are

$$dom(x_{av_1}) = \{any\ value\ between\ 0\ and\ 5\}$$
$$dom(x_{av_2}) = \{any\ value\ between\ 0\ and\ 5\}$$
$$dom(x_{av_3}) = \{any\ value\ between\ 0\ and\ 5\}$$

By solving this problem as a CSP, the core of the game is

$$Core = \{(x_{av_1}, x_{av_2}, x_{av_3}) : \sum_{i=1}^{i=3} x_{av_i} = f(N),\ and$$
$$x(S) \geq f(S),\ \forall\ S \in N\}$$

$$Core = \{(\$2, \$0, \$3)(\$3, \$0, \$2)(\$4, \$0, \$1)(\$5, \$0, \$0)\}.$$

## 5 Simulation Tests

Simulation was conducted to test and validate the performance of the CSOF. The simulation was performed in MATLAB using the concept of Object Oriented Programming (OOP). A two-lane roadway sub-network containing three SIs was chosen in Waterloo, ON, Canada, to conduct the simulation (Fig. 3). The SIs are as follows: SI1, Westmount Road North with Columbia Street West; SI2, Westmount Road North with Bearinger Road; and SI3, Northfield Drive West with Weber Street North. Every SI has a static TL system such that TL1, TL2, and TL3 for SI1, SI2, and SI3 respectively. Each TL control system has a two-phase cycle where the East-West roadways are one phase and South-North roadways are the other phase. Each phase has a signal design of Green-Yellow-Red; however, for simplicity, the yellow-light time is assumed to be part of the green-light-time duration. To enhance safety, one second of red-light time is given to all the roadways between every two consecutive phases.

In order to overcome randomization and capture the real behaviour of traffic, the simulation was run for more than three hours. The maximum and minimum speed limits on any roadway in the network are $v_{max} = 60 km/hour$ and $v_{min} = 10 km/hour$ respectively where the highest volume of traffic a road segment may have is assumed to be 85percent of the maximum density, $D_{max}(L_j)$. AVs are generated randomly into the network based on Poisson Distribution. The generated AVs travel at an average speed of $50 km/hour$ until they get within the activation distance (i.e., the distance at which the vehicles get within the V2I communication range and start cooperating). The activation distance was fixed at $500 meters$. The performance of the CSOF is compared to a Non-Cooperative Speed Optimization algorithm (NCSO) (i.e., the vehicles individually compute their optimal speeds). Once they are within range, AVs start conducting speed optimization based on the functioning optimization technique (i.e., CSOF or NCSO) to have a chance of meeting the green light when they arrive at the TL.

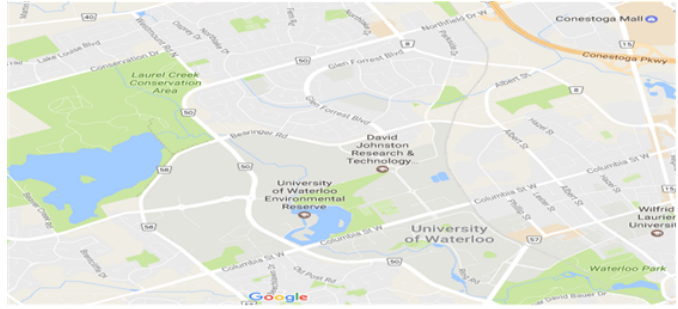

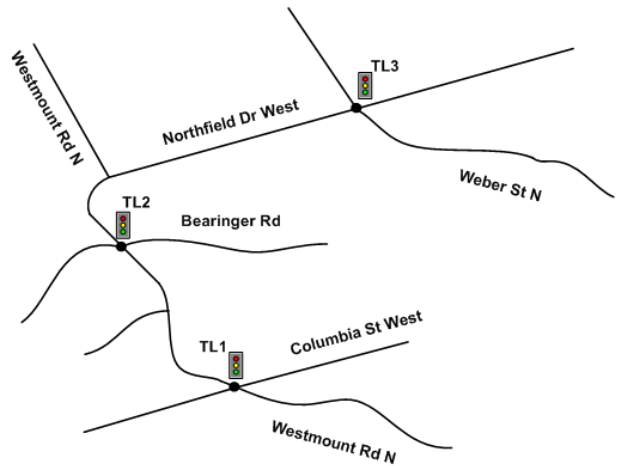

**Figure 3:** A sub-network with three signalized intersections.

Figs. 4, 5, and 6 report the results of the total average idling time, total average number of stops, and total average energy consumption at SI1, SI2, and SI3, comparing the CSOF to the NCSO algorithm. As can be seen, the CSOF outperformed the NCSO by achieving lower average idling times and average number of stops. This is because the conflicting passing times of vehicles through the intersections are resolved by the CSOF. All the AVs using the CSOF, meant to arrive during the green-light time, are al-

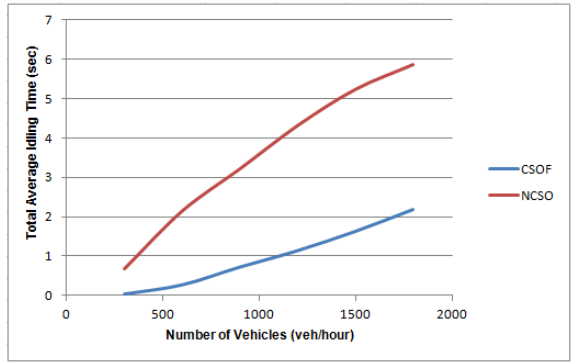

**Figure 4:** Total average idling time at SI1, SI2, and SI3.

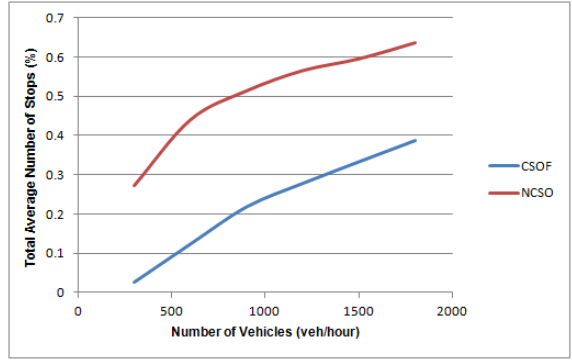

**Figure 5:** Total average number of stops at SI1, SI2, and SI3.

located time tokens before reaching the intersections. As such, when they arrive, they are able to pass through smoothly during their allocated times. In addition, due to the road and signal timing constraints, the AVs that could only arrive at the intersections during the red-light times were not allocated time tokens before reaching the intersections. These AVs joined the queues with less waiting times. As mentioned previously, the time needed to clear the queue is excluded from that available as time tokens to the approaching AVs. It can be noticed from the result figures that for the two compared techniques, as the number of AVs approaching the intersections increases, the average idling times and number of stops become greater, and the improvement achieved by the CSOF in minimizing the average idling times is more significant. The reductions in average idling times that have been achieved by the CSOF when compared to the NCSO for SI1, SI2, and SI3 are summarized in Tables 3, 4, and 5 respectively.

In addition, the CSOF has achieved lower average energy consumption. In the case of CSOF, as soon as a vehicle is allocated a time token, it maintains its speed so that it passes the intersection smoothly during its allocated token. Hence, in general there are less speed variations resulting from the AVs using the CSOF. To prove this, Fig. 7 captures the speed trajectories of six AVs approaching SI1 from the Westmount Road North direction from the moment they joined the activation distance of the CSOF until they passed the stop line of the intersection. It is clear that on average, the reductions in speed variations of AVs using the CSOF are not significant as compared to those using the NCSO algorithm. As a result, the energy savings of AVs using the CSOF have not been significant as reported in Fig. 6.

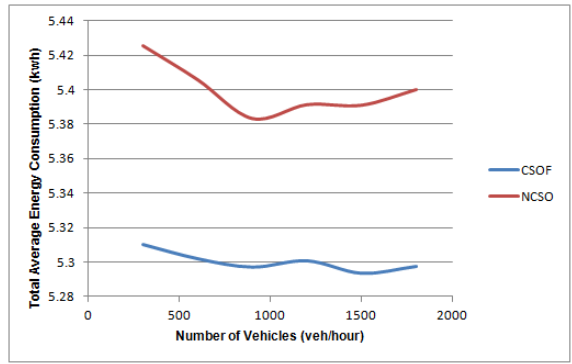

**Figure 6:** Total average energy consumption of vehicles approaching SI1, SI2, and SI3.

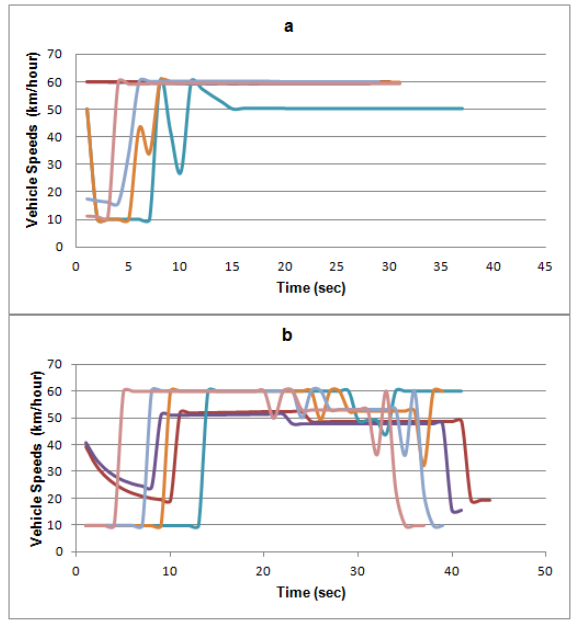

**Figure 7:** Speed trajectories of six vehicles approaching intersection 1, (a) vehicles using the CSOF, (b) vehicles using the NCSO.

Furthermore, the total average idling time and number of stops achieved by the

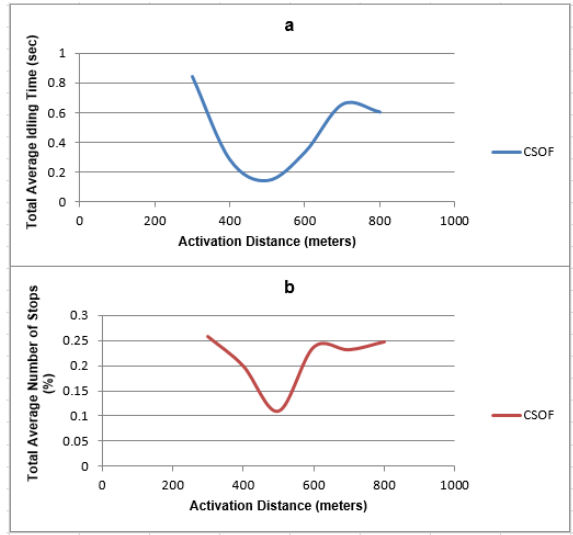

**Figure 8:** Activation distance analysis for intersection 2, (a) total average idling time, (b) total average number of stops.

CSOF were investigated with respect to the V2I activation distance. SI2 was chosen to conduct the investigation. It is assumed that the V2I communication radio is available in a range of up to 800*meters* away from the intersection, so it was varied from 300*meters* to 800*meters* in steps of 100*meters*. Fig. 8 depicts the total values of average idling time and average number of stops being achieved by the CSOF. As can be seen, it is concluded that the optimal point of activation is found near 500*meters*. At shorter distances, some AVs are forced to arrive during the red-light time due to the fact that the time available to the AVs to get allocated tokens, or reallocated tokens after playing a game, to adjust their speeds accordingly is not enough to result in low average values of idling time and number of stops. At further activation distances, the average values of idling time and number of stops are slightly increased but remain near the optimal level.

**Table 3:** Reduction in average idling time at signalized intersection 1.

| AVs (veh/hour) | 300 | 600 | 900 | 1200 | 1500 | 1800 |
|---|---|---|---|---|---|---|
| NCSO (sec) | 0.829 | 3.831 | 4.301 | 5.55 | 5.967 | 6.262 |
| CSOF (sec) | 0.05 | 0.441 | 1.227 | 1.688 | 2.521 | 2.933 |
| Reduction (%) | 94 | 88 | 71 | 70 | 58 | 53 |

**Table 4:** Reduction in average idling time at signalized intersection 2.

| AVs (veh/hour) | 300 | 600 | 900 | 1200 | 1500 | 1800 |
|---|---|---|---|---|---|---|
| NCSO (sec) | 0.388 | 0.888 | 2.299 | 2.754 | 3.614 | 4.632 |
| CSOF (sec) | 0.035 | 0.176 | 0.347 | 0.603 | 0.725 | 1.08 |
| Reduction (%) | 91 | 80 | 85 | 78 | 80 | 77 |

**Table 5:** Reduction in average idling time at signalized intersection 3.

| AVs (veh/hour) | 300 | 600 | 900 | 1200 | 1500 | 1800 |
|---|---|---|---|---|---|---|
| NCSO (sec) | 0.711 | 1.826 | 2.694 | 4.649 | 6.208 | 6.412 |
| CSOF (sec) | 0.044 | 0.159 | 0.483 | 1.084 | 1.727 | 2.878 |
| Reduction (%) | 94 | 91 | 82 | 77 | 72 | 55 |

# 6 Conclusion

This research has addressed the problem of minimizing the total average idling times and number of stops for Electric Autonomous Vehicles (EAVs) approaching Signalized Intersections (SIs). A Cooperative Speed Optimization Framework (CSOF) was proposed. The proposed framework consists of three modules to tackle issues of individual speed optimization, conflict recognition, and cooperative speed decision making. Simulation was conducted to investigate the performance of the CSOF in comparison with a Non-Cooperative Speed Optimization (NCSO) algorithm. The simulation tested the performance of the two techniques under various traffic conditions, showing that the CSOF outperformed the NCSO in terms of minimizing the total average idling times, number of stops, and energy consumption. Furthermore, the performance of the CSOF was investigated in terms of the V2I communication range. It was concluded that the EAVs using the CSOF could achieve the minimum values of total average idling times and number of stops at a communication range activation near 500*meters*. The future work of this research should investigate the scalability of the CSOF under various geometrical designs of SIs and TL phase/cycle settings. Moreover, the cooperation of a dynamic traffic light system with the EAVs using the CSOF should be investigated to achieve further minimization of average idling times, number of stops, and energy consumption.